\documentclass{article}

\begin{document}    
    
\input{amssym.def}    
\input{amssym.tex}    

\title{Noncommutative geometry of algebraic curves}

\author{Igor  ~Nikolaev
\footnote{Partially supported 
by NSERC.}}

\date{}    
 \maketitle    
    
\newtheorem{thm}{Theorem}    
\newtheorem{lem}{Lemma}    
\newtheorem{dfn}{Definition}    
\newtheorem{rmk}{Remark}    
\newtheorem{cor}{Corollary}    
\newtheorem{prp}{Proposition}    
\newtheorem{cnj}{Conjecture}   
\newtheorem{prb}{Problem}

 \begin{abstract}    
A covariant functor from the category of generic complex algebraic curves to
a category of the $AF$-algebras is constructed. The construction is
based on a representation of the Teichm\"uller space of a curve
by the measured foliations due to Douady, Hubbard, Masur and Thurston.
The functor maps isomorphic algebraic curves
to the stably isomorphic $AF$-algebras.

\vspace{7mm}    
    

\vspace{5mm}    
{\it AMS (MOS) Subj. Class.:  14H10, 46L40, 58F10}    
\end{abstract}

\section*{Introduction}    
There is an  interest in the last decade in a dictionary
between commutative geometry of the rings of polynomials in the complex projective plane  
and noncommutative geometry of the rings of operators on a   Hilbert space. 
A  special  case of  elliptic curves and noncommutative tori suggests
that such a  dictionary exists \cite{Kon1}, \cite{Man1},  \cite{PoSch1}, \cite{SoVo1} 
{\it et al.}  In the present note we construct a functor
from the  category of generic  complex algebraic curves to a category 
of the operator algebras, known as the $AF$-algebras, which realizes such
a dictionary. The functor maps  isomorphic algebraic curves to the stably
isomorphic (Morita equivalent) $AF$-algebras. This fact, interesting
on its own, has applications, e.g. in the construction of a faithful 
representation of the mapping class group.

Let us outline our construction.  Denote by $T_S(g)$ the Teichm\"uller
space of genus $g\ge 1$ with a distinguished point $S$. Let
$q\in H^0(S, \Omega^{\otimes 2})$ be a holomorphic quadratic 
differential on the Riemann surface $S$, such that all zeroes
of $q$ (if any) are simple. By $\widetilde S$ we understand a double 
cover of $S$ ramified over the zeroes of $q$ and by
$H_1^{odd}(\widetilde S)$ the odd part of the integral homology of $\widetilde S$ relatively the  zeroes.
Note that $H_1^{odd}(\widetilde S)\cong {\Bbb Z}^{n}$, where $n=6g-6$ if $g\ge 2$ and $n=2$ if $g=1$. 
The fundamental result of Hubbard and Masur \cite{HuMa1} says that
$
T_S(g)\cong Hom~(H_1^{odd}(\widetilde S); {\Bbb R})-\{0\},
$
where $0$ is the zero homomorphism.

Finally, denote by $\lambda=(\lambda_1,\dots,\lambda_{n})$ the image of a basis of 
$H_1^{odd}(\widetilde S)$ in the real line ${\Bbb R}$, such that $\lambda_1\ne0$. 
(Note that such an option always exists, since the zero homomorphism is excluded.) 
We let $\theta=(\theta_1,\dots,\theta_{n-1})$, where $\theta_i=\lambda_{i-1}/\lambda_1$. 
Recall that,  up to a scalar multiple, vector $(1,\theta)\in {\Bbb R}^{n}$ is the limit
of a generically convergent Jacobi-Perron continued fraction \cite{B}:
\begin{equation}\label{eq2}
\left(\matrix{1\cr \theta}\right)=
\lim_{k\to\infty} \left(\matrix{0 & 1\cr I & b_1}\right)\dots
\left(\matrix{0 & 1\cr I & b_k}\right)
\left(\matrix{0\cr {\Bbb I}}\right),
\end{equation}
where $b_i=(b^{(i)}_1,\dots, b^{(i)}_{n-1})^T$ is a vector of the non-negative integers,  
$I$ the unit matrix and ${\Bbb I}=(0,\dots, 0, 1)^T$. We introduce an  $AF$-algebra \cite{E}, 
${\Bbb A}_{\theta}$, via the Bratteli diagram:

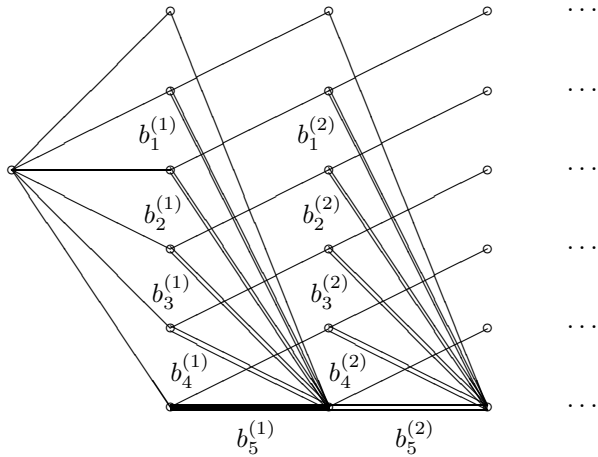
\begin{figure}[here]
\begin{picture}(300,200)(0,0)

\put(40,110){\circle{3}}

\put(40,110){\line(1,1){60}}
\put(40,110){\line(2,1){60}}
\put(40,110){\line(1,0){60}}
\put(40,110){\line(2,-1){60}}
\put(40,110){\line(1,-1){60}}
\put(40,110){\line(2,-3){60}}

\put(100,20){\circle{3}}
\put(100,50){\circle{3}}
\put(100,80){\circle{3}}
\put(100,110){\circle{3}}
\put(100,140){\circle{3}}
\put(100,170){\circle{3}}

\put(160,20){\circle{3}}
\put(160,50){\circle{3}}
\put(160,80){\circle{3}}
\put(160,110){\circle{3}}
\put(160,140){\circle{3}}
\put(160,170){\circle{3}}

\put(220,20){\circle{3}}
\put(220,50){\circle{3}}
\put(220,80){\circle{3}}
\put(220,110){\circle{3}}
\put(220,140){\circle{3}}
\put(220,170){\circle{3}}

\put(160,20){\line(2,1){60}}
\put(160,19){\line(1,0){60}}
\put(160,21){\line(1,0){60}}
\put(160,50){\line(2,1){60}}
\put(160,49){\line(2,-1){60}}
\put(160,51){\line(2,-1){60}}
\put(160,80){\line(2,1){60}}
\put(160,79){\line(1,-1){60}}
\put(160,81){\line(1,-1){60}}
\put(160,110){\line(2,1){60}}
\put(160,109){\line(2,-3){60}}
\put(160,111){\line(2,-3){60}}
\put(160,140){\line(2,1){60}}
\put(160,139){\line(1,-2){60}}
\put(160,141){\line(1,-2){60}}
\put(160,170){\line(2,-5){60}}


\put(100,20){\line(2,1){60}}
\put(100,19){\line(1,0){60}}
\put(100,20){\line(1,0){60}}
\put(100,21){\line(1,0){60}}
\put(100,50){\line(2,1){60}}
\put(100,49){\line(2,-1){60}}
\put(100,51){\line(2,-1){60}}
\put(100,80){\line(2,1){60}}
\put(100,79){\line(1,-1){60}}
\put(100,81){\line(1,-1){60}}
\put(100,110){\line(2,1){60}}
\put(100,109){\line(2,-3){60}}
\put(100,111){\line(2,-3){60}}
\put(100,140){\line(2,1){60}}
\put(100,139){\line(1,-2){60}}
\put(100,141){\line(1,-2){60}}
\put(100,170){\line(2,-5){60}}


\put(250,20){$\dots$}
\put(250,50){$\dots$}
\put(250,80){$\dots$}
\put(250,110){$\dots$}
\put(250,140){$\dots$}
\put(250,170){$\dots$}


\put(125,5){$b_5^{(1)}$}
\put(100,30){$b_4^{(1)}$}
\put(93,60){$b_3^{(1)}$}
\put(90,90){$b_2^{(1)}$}
\put(88,120){$b_1^{(1)}$}

\put(185,5){$b_5^{(2)}$}
\put(160,30){$b_4^{(2)}$}
\put(153,60){$b_3^{(2)}$}
\put(150,90){$b_2^{(2)}$}
\put(148,120){$b_1^{(2)}$}


\end{picture}

\caption{The Bratteli diagram of a toric $AF$-algebra ${\Bbb A}_{\theta}$ (case  $g=2$).}
\end{figure}

\noindent
where numbers $b_j^{(i)}$  indicate the multiplicity  of edges of 
the graph. Let us call ${\Bbb A}_{\theta}$ a {\it toric $AF$-algebra}. 
Note that in the case $g=1$, the Jacobi-Perron 
fraction coincides with the  regular continued fraction and ${\Bbb A}_{\theta}$
becomes the Effros-Shen $AF$-algebra of a noncommutative torus
\cite{EfSh1}.  Roughly, the question addressed in this note
can be formulated as follows.

\medskip\noindent
{\bf Main problem.} {\it Suppose that $C, C'\in T_S(g)$ are isomoprhic complex algebraic 
curves of genus $g\ge 1$. Find  an equivalence relation between  the corresponding toric $AF$-algebras 
${\Bbb A}_{\theta}, {\Bbb A}_{\theta'}$.}

\medskip\noindent
It is remarkable that the  difficult  question  of an isomorphism   between the two 
algebraic  curves  has an amazingly  simple answer in terms of the operator algebras. Recall that the fundamental
equivalence in noncommutative geometry is a stable  isomorphism
of the operator algebras rather than an  isomorphism.  The operator algebras
${\Bbb A}, {\Bbb A}'$ are stably isomorphic whenever ${\Bbb A}\otimes {\cal K}$
is isomorphic to  ${\Bbb A}'\otimes {\cal K}$, where ${\cal K}$ is the $C^*$-algebra
of compact operators. A short  answer to the main problem will be that ${\Bbb A}_{\theta}, {\Bbb A}_{\theta'}$
are stably isomorphic. Unfortunately, this beautiful fact holds only for the  typical (generic)
algebraic curves, which we will further specify.  An obstacle in the remaining cases  
is a (non-generic) divergence of the Jacobi-Perron algorithm. One always has  
the  convergence only in the case $g=1$, i.e. when ${\Bbb A}_{\theta}$ is an Effros-Shen algebra.

Denote by $V$ the maximal subset of $T_S(g)$ such that for every complex curve $C\in V$, 
there exists a convergent Jacobi-Perron continued fraction. 
Let $F$ be the map which sends the complex algebraic curves into the toric $AF$-algebras
according to the formula  $C\mapsto {\Bbb A}_{\theta}$ described  above.  
Finally, let $W$ be the image of  $V$ under the mapping $F$.
A  summary of  our results is contained  in the following theorem.
\begin{thm}\label{thm1}
The set $V$ is a generic subset of $T_S(g)$ and the  map $F$ has the following properties:

\smallskip
(i) $V\cong W\times (0,\infty)$ is a trivial fiber bundle, whose
projection  map $p: V\to W$ coincides with $F$;

\smallskip
(ii) $F: V\to W$ is a covariant functor, which maps isomorphic complex algebraic curves 
$C,C'\in V$ to the stably isomorphic toric $AF$-algebras ${\Bbb A}_{\theta},{\Bbb A}_{\theta'}\in W$.
\end{thm}

\medskip\noindent    
The  note is organized as follows.  In section 1 some useful definitions are introduced and 
theorem \ref{thm1}  is proved in section 2.

\section{Preliminaries}
This section is a brief review of the results necessary to prove
our main theorem. For a systematic account,  we refer the reader
to \cite{B}, \cite{HuMa1} and \cite{Thu1}.

\subsection{Measured foliations and $T(g)$}
{\bf A.}  A measured foliation, ${\cal F}$, on a surface $X$
is a  partition of $X$ into the singular points $x_1,\dots,x_n$ of
order $k_1,\dots, k_n$ and regular leaves (1-dimensional submanifolds). 
On each  open cover $U_i$ of $X-\{x_1,\dots,x_n\}$ there exists a non-vanishing
real-valued closed 1-form $\phi_i$  such that: 
(i)  $\phi_i=\pm \phi_j$ on $U_i\cap U_j$;
(ii) at each $x_i$ there exists a local chart $(u,v):V\to {\Bbb R}^2$
such that for $z=u+iv$, it holds $\phi_i=Im~(z^{k_i\over 2}dz)$ on
$V\cap U_i$ for some branch of $z^{k_i\over 2}$. 
The pair $(U_i,\phi_i)$ is called an atlas for measured foliation ${\cal F}$.
Finally, a measure $\mu$ is assigned to each segment $(t_0,t)\in U_i$, which is  transverse to
the leaves of ${\cal F}$, via the integral $\mu(t_0,t)=\int_{t_0}^t\phi_i$. The 
measure is invariant along the leaves of ${\cal F}$, hence the name.

\smallskip\noindent
{\bf B.} Let $S$ be a Riemann surface, and $q\in H^0(S,\Omega^{\otimes 2})$ a holomorphic
quadratic differential on $S$. The lines $Re~q=0$ and $Im~q=0$ define a pair
of measured foliations on $R$, which are transversal to each other outside the set of 
singular points. The set of singular points is common to both foliations and coincides
with the zeroes of $q$. The above measured foliations are said to represent the  vertical and horizontal 
 trajectory structure  of $q$, respectively.

\smallskip\noindent
{\bf C.}  Let $T(g)$ be the Teichm\"uller space of the topological surface $X$ of genus $g\ge 1$,
i.e. the space of the complex structures on $X$. 
Consider the vector bundle $p: Q\to T(g)$ over $T(g)$ whose fiber above a point 
$S\in T_g$ is the vector space $H^0(S,\Omega^{\otimes 2})$.   
Given non-zero $q\in Q$ above $S$, we can consider horizontal measured foliation
${\cal F}_q\in \Phi_X$ of $q$, where $\Phi_X$ denotes the space of equivalence
classes of measured foliations on $X$. If $\{0\}$ is the zero section of $Q$,
the above construction defines a map $Q-\{0\}\longrightarrow \Phi_X$. 
For any ${\cal F}\in\Phi_X$, let $E_{\cal F}\subset Q-\{0\}$ be the fiber
above ${\cal F}$. In other words, $E_{\cal F}$ is a subspace of the holomorphic 
quadratic forms whose horizontal trajectory structure coincides with the 
measured foliation ${\cal F}$. 
Note that, if ${\cal F}$ is a measured foliation with the simple zeroes (a generic case),  
then $E_{\cal F}\cong {\Bbb R}^n - 0$, while $T(g)\cong {\Bbb R}^n$, where $n=6g-6$ if
$g\ge 2$ and $n=2$ if $g=1$.

\bigskip\noindent
{\bf Theorem (Hubbard-Masur \cite{HuMa1})}
{\it The restriction of $p$ to $E_{\cal F}$ defines a homeomorphism (an embedding)
$h_{\cal F}: E_{\cal F}\to T(g)$.}

\bigskip\noindent
{\bf D.}
The Hubbard-Masur result implies that the measured foliations  parametrize  
the space $T(g)-\{pt\}$, where $pt= h_{\cal F}(0)$.
Indeed, denote by  ${\cal F}'$ a vertical trajectory structure of  $q$. Since ${\cal F}$
and ${\cal F}'$ define $q$, and ${\cal F}=Const$ for all $q\in E_{\cal F}$, one gets a homeomorphism 
between $T(g)-\{pt\}$ and  $\Phi_X$, where $\Phi_X\cong {\Bbb R}^n-0$ is the space of 
equivalence classes of the measured foliations ${\cal F}'$ on $X$. 
Note that the above parametrization depends on a foliation ${\cal F}$.
However, there exists a unique canonical homeomorphism $h=h_{\cal F}$
as follows. Let $Sp ~(S)$ be the length spectrum of the Riemann surface
$S$ and $Sp ~({\cal F}')$ be the set positive reals $\inf \mu(\gamma_i)$,
where $\gamma_i$ runs over all simple closed curves, which are transverse 
to the foliation ${\cal F}'$. A canonical homeomorphism 
$h=h_{\cal F}: \Phi_X\to T(g)-\{pt\}$ is defined by the formula
$Sp ~({\cal F}')= Sp ~(h_{\cal F}({\cal F}'))$ for $\forall {\cal F}'\in\Phi_X$. 
Thus, the following corollary is true.
\begin{cor}\label{cor1}
There exists a unique homeomorphism $h:\Phi_X\to T(g)-\{pt\}$.
\end{cor}

\medskip\noindent
{\bf E.} 
Recall that $\Phi_X$ is the space of equivalence classes of measured
foliations on the topological surface $X$. Following Douady and Hubbard
\cite{DoHu1}, we consider a coordinate system on $\Phi_X$,
suitable for the proof of theorem \ref{thm1}. 
For clarity, let us make a generic assumption that $q\in H^0(S,\Omega^{\otimes 2})$
is a non-trivial holomorphic quadratic differential with only simple zeroes. 
We wish to construct a Riemann surface of $\sqrt{q}$, which is a double cover
of $S$ with ramification over the zeroes of $q$. Such a surface, denoted by
$\widetilde S$, is unique and has an advantage of carrying a holomorphic
differential $\omega$, such that $\omega^2=q$. We further denote by 
$\pi:\widetilde S\to S$ the covering projection. The vector space
$H^0(\widetilde S,\Omega)$ splits into the direct sum
$H^0_{even}(\widetilde S,\Omega)\oplus H^0_{odd}(\widetilde S,\Omega)$
in view of  the involution $\pi^{-1}$ of $\widetilde S$, and
the vector space $H^0(S,\Omega^{\otimes 2})\cong H^0_{odd}(\widetilde S,\Omega)$.
Let $H_1^{odd}(\widetilde S)$ be the odd part of the homology of $\widetilde S$
relatively  the zeroes of $q$.   Consider the pairing
$H_1^{odd}(\widetilde S)\times H^0(S, \Omega^{\otimes 2})\to {\Bbb C}$,
defined by the integration  $(\gamma, q)\mapsto \int_{\gamma}\omega$. 
We shall take the associated map
$\psi_q: H^0(S,\Omega^{\otimes 2})\to Hom~(H_1^{odd}(\widetilde S); {\Bbb C})$
and let $h_q= Re~\psi_q$. 
\begin{lem}{\bf (Douady-Hubbard \cite{DoHu1})}
The map
$
h_q: H^0(S, \Omega^{\otimes 2})\longrightarrow Hom~(H_1^{odd}(\widetilde S); {\Bbb R})
$
is an ${\Bbb R}$-isomorphism. 
\end{lem}
Since  each  ${\cal F}\in \Phi_X$ is the  vertical foliation 
$Re~q=0$ for a $q\in H^0(S, \Omega^{\otimes 2})$, the Douady-Hubbard lemma
implies that $\Phi_X\cong Hom~(H_1^{odd}(\widetilde S); {\Bbb R})$.
By  formulas for the relative homology, 
one finds that $H_1^{odd}(\widetilde S)\cong {\Bbb Z}^{n}$,
where $n=6g-6$ if $g\ge 2$ and $n=2$ if $g=1$.     
Finally, each $h\in Hom~({\Bbb Z}^{n}; {\Bbb R})$ is given
by the reals  $\lambda_1=h(e_1),\dots, \lambda_{n}=h(e_{n})$,
where $(e_1,\dots, e_{n})$ is a basis in ${\Bbb Z}^{n}$.
The numbers   $(\lambda_1,\dots,\lambda_{n})$ are the coordinates in the space $\Phi_X$
and, in view of the corollary \ref{cor1}, in  the Teichm\"uller space $T(g)$.

\subsection{The Jacobi-Perron continued fraction}
{\bf A.}
Let $a_1,a_2\in {\Bbb N}$ such that $a_2\le a_1$. Recall that the greatest common
divisor of $a_1,a_2$, $GCD(a_1,a_2)$, can be determined from the Euclidean algorithm:
\begin{equation}\label{eq5}
\left\{
\begin{array}{cc}
a_1 &= a_2b_1 +r_3\nonumber\\
a_2 &= r_3b_2 +r_4\nonumber\\
r_3 &= r_4b_3 +r_5\nonumber\\
\vdots & \nonumber\\
r_{k-3} &= r_{k-2}b_{k-1}+r_{k-1}\nonumber\\
r_{k-2} &= r_{k-1}b_k,
\end{array}
\right.
\end{equation}
where $b_i\in {\Bbb N}$ and $GCD(a_1,a_2)=r_{k-1}$. 
The Euclidean algorithm can be written as the regular continued 
fraction
\begin{equation}\label{eq6}
\theta={a_1\over a_2}=b_1+{1\over\displaystyle b_2+
{\strut 1\over\displaystyle +\dots+ {1\over b_k}}}
=[b_1,\dots b_k].
\end{equation}
If $a_1, a_2$ are non-commensurable, in the sense that $\theta\in {\Bbb R}-{\Bbb Q}$,
then the Euclidean algorithm never stops and $\theta=[b_1, b_2, \dots]$. Note that the regular  
continued fraction can be written in the matrix form:
\begin{equation}\label{eq7}
\left(\matrix{1\cr \theta}\right)=
\lim_{k\to\infty} \left(\matrix{0 & 1\cr 1 & b_1}\right)\dots
\left(\matrix{0 & 1\cr 1 & b_k}\right)
\left(\matrix{0\cr 1}\right). 
\end{equation}

\medskip\noindent
{\bf B.}
The Jacobi-Perron algorithm and connected (multidimensional) continued 
fraction generalizes the Euclidean algorithm to the case $GCD(a_1,\dots,a_n)$
when $n\ge 2$. Namely, let $\lambda=(\lambda_1,\dots,\lambda_n)$,
$\lambda_i\in {\Bbb R}-{\Bbb Q}$ and  $\theta_{i-1}={\lambda_i\over\lambda_1}$, where
$1\le i\le n$.   The continued fraction 
$$
\left(\matrix{1\cr \theta_1\cr\vdots\cr\theta_{n-1}} \right)=
\lim_{k\to\infty} 
\left(\matrix{0 &  0 & \dots & 0 & 1\cr
              1 &  0 & \dots & 0 & b_1^{(1)}\cr
              \vdots &\vdots & &\vdots &\vdots\cr
              0 &  0 & \dots & 1 & b_{n-1}^{(1)}}\right)
\dots 
\left(\matrix{0 &  0 & \dots & 0 & 1\cr
              1 &  0 & \dots & 0 & b_1^{(k)}\cr
              \vdots &\vdots & &\vdots &\vdots\cr
              0 &  0 & \dots & 1 & b_{n-1}^{(k)}}\right)
\left(\matrix{0\cr 0\cr\vdots\cr 1} \right),
$$
where $b_i^{(j)}\in {\Bbb N}\cup\{0\}$, is called the {\it Jacobi-Perron
algorithm (JPA)}. Unlike the regular continued fraction algorithm,
the JPA may diverge for certain vectors $\lambda\in {\Bbb R}^n$. However, 
for points of a generic subset of ${\Bbb R}^n$, the JPA converges. 
We  characterize such a set in the next paragraph.

\medskip\noindent
{\bf C.}
The convergence of the JPA algorithm can be characterized in terms of
the measured foliations. Let ${\cal F}\in\Phi_X$ be a measured foliation
on the surface $X$ of genus $g\ge 1$. Recall that ${\cal F}$ is called uniquely
ergodic if every invariant measure of ${\cal F}$ is a multiple
of the Lebesgue measure. By the Masur-Veech theorem, there exists
a generic subset $V\subset \Phi_X$ such that each ${\cal F}\in V$
is uniquely ergodic \cite{Mas1}, \cite{Vee1}.
We let $\lambda=(\lambda_1,\dots,\lambda_{n})$ be the vector with
coordinates $\lambda_i=\mu ({\gamma_i})$, where $\gamma_i\in H_1^{odd}(\widetilde S)$,
see \S 1.1.E. By an abuse of notation, we shall say that $\lambda\in V$. 
By a duality  between the measured foliations and the interval exchange transformations \cite{Mas1}, 
the following characterization of convergence of the JPA is true.
\begin{lem}
{\bf (Bauer \cite{Bau1})}
The JPA converges if and only if $\lambda\in V\subset {\Bbb R}^{n}$.  
\end{lem}

\section{Proof of theorem 1}
Let us outline the proof. We consider the following categories:
(i) generic complex algebraic curves ${\cal V}=V\subset T(g)$;
(ii) pseudo-lattices ${\cal PL}$; (iii) projective pseudo-lattices
${\cal PPL}$ and (iv)  category ${\cal W}$ of the toric $AF$-algebras.
First, we show that ${\cal V}\cong {\cal PL}$ are equivalent categories,  such that isomorphic
complex algebraic  curves $C,C'\in {\cal V}$ map to isomorphic pseudo-lattices
$PL, PL'\in {\cal PL}$.  Next, a non-injective functor $F: {\cal PL}\to {\cal PPL}$
is constructed.  The $F$ maps isomorphic pseudo-lattices to isomorphic
projective pseudo-lattices and $Ker ~F\cong (0,\infty)$. 
Finally, it is shown that a subcategory $U\subseteq {\cal PPL}$ and ${\cal W}$ are the 
equivalent categories. In other words, we have the following diagram:
\begin{equation}\label{eq8}
{\cal V}
\buildrel\rm\alpha
\over\longrightarrow
{\cal PL}
\buildrel\rm F
\over\longrightarrow
U
\buildrel\rm \beta
\over\longrightarrow
{\cal W},
\end{equation}
where $\alpha$ is an injective map, $\beta$ is a bijection  and $Ker ~F\cong (0,\infty)$.

\bigskip\noindent
\underline{Category ${\cal V}$}. 
Let $Mod ~X$ be the mapping class group of the surface $X$. 
A {\it complex algebraic curve}  is a triple $(X, C, j)$, where
$X$ is a topological surface of genus $g\ge 1$,  $j: X\to C$ is a
complex (conformal) parametrization of $X$ and $C$ is a Riemann surface.
A {\it morphism} of complex algebraic curves  
$(X, C, j)\to (X, C', j')$   is the identity
$j\circ\psi=\varphi\circ j'$, 
where $\varphi\in Mod~X$ is a diffeomorphism of $X$  and $\psi$ is an isomorphism of Riemann surfaces. 
A  category of generic complex algebraic curves, ${\cal V}$, consists of $Ob~({\cal V})$
which are complex algebraic curves  $C\in V\subset T_S(g)$ and morphisms $H(C,C')$
between $C,C'\in Ob~({\cal V})$ which coincide with the morphisms
specified above. For any  $C,C',C''\in Ob~({\cal V})$
and any morphisms $\varphi': C\to C'$, $\varphi'': C'\to C''$ a 
morphism $\phi: C\to C''$ is the composite of $\varphi'$ and
$\varphi''$, which we write as $\phi=\varphi''\varphi'$. 
The identity  morphism, $1_C$, is a morphism $H(C,C)$.

\bigskip\noindent
\underline{Category ${\cal PL}$}. 
A {\it pseudo-lattice}
\footnote{See \cite{Man1}.}
 (of rank $n$) is a triple $(\Lambda, {\Bbb R}, j)$, where
$\Lambda\cong {\Bbb Z}^n$ and $j: \Lambda\to {\Bbb R}$ is a homomorphism.  
A morphism of pseudo-lattices $(\Lambda, {\Bbb R}, j)\to (\Lambda, {\Bbb R}, j')$
is the identity
$j\circ\psi=\varphi\circ j'$,
where $\varphi$ is a group homomorphism and $\psi$ is an inclusion 
map, i.e. $j'(\Lambda')\subseteq j(\Lambda)$.  
Any isomorphism class of a pseudo-lattice contains
a representative given by $j: {\Bbb Z}^n\to  {\Bbb R}$ such that 
$j(1,0,\dots, 0)=\lambda_1, \quad j(0,1,\dots,0)=\lambda_2,\quad \dots, \quad  j(0,0,\dots,1)=\lambda_n,$
where $\lambda_1,\lambda_2,\dots,\lambda_n$ are positive reals.
The pseudo-lattices of rank $n$  make up a category, which we denote by ${\cal PL}_n$. 
\begin{lem}{\bf (basic lemma)}\label{lm5}
Let $g\ge 2$ ($g=1$) and $n=6g-6$ ($n=2$). 
There exists an injective covariant functor $\alpha: {\cal V}\to {\cal PL}_{n}$,
which maps isomorphic complex algebraic curves $C,C'\in {\cal V}$ to the isomorphic
pseudo-lattices $PL,PL'\in {\cal PL}_{n}$.
\end{lem}
{\it Proof.} Let $\alpha: T(g)-\{pt\}\to Hom~(H_1^{odd}(\widetilde S); {\Bbb R})-\{0\}$
be the map  constructed in \S 1.1. By the corollary from the Hubbard-Masur theorem, $\alpha$
is a homeomorphism. In particular, $\alpha$ is an injective map.

Let us find the image $\alpha(\varphi)\in Mor~({\cal PL})$ of 
$\varphi\in Mor~({\cal V})$. Let $\varphi\in Mod~X$ be a diffeomorphism
of $X$, and let $\widetilde X\to X$ be the ramified double cover of $X$
as explained in \S 1.1.E. We denote by $\widetilde\varphi$ the induced
mapping on $\widetilde X$. Note that $\widetilde\varphi$ is a diffeomorphism
of $\widetilde X$ modulo the covering involution ${\Bbb Z}_2$. Denote by
$\widetilde\varphi^*$ the action of $\widetilde\varphi$ on 
$H_1^{odd}(\widetilde X)\cong {\Bbb Z}^{n}$. Since $\widetilde\varphi~mod~{\Bbb Z}_2$ is
a diffeomorphism of $\widetilde X$, $\widetilde\varphi^*\in GL_{n}({\Bbb Z})$. 
Thus, $\alpha(\varphi)=\widetilde\varphi^*\in Mor~({\cal PL})$.

Let us show that $\alpha$ is a functor.  Indeed, let $C,C'\in {\cal V}$ be the isomorphic
complex algebraic curves, such that $C'=\varphi(C)$ for a $\varphi\in Mod~X$. 
Let $a_{ij}$ be the elements of matrix $\widetilde\varphi^*\in GL_{n}({\Bbb Z})$. 
Recall that $\lambda_i=\int_{\gamma_i}\phi$ for a closed 1-form $\phi= Re~\omega$
and $\gamma_i\in H_1^{odd}(\widetilde X)$. Then
$
\gamma_j=\sum_{i=1}^{n} a_{ij}\gamma_i, \quad j=1,\dots, n 
$
are the elements of a new basis in $H_1^{odd}(\widetilde X)$. By the integration rules,
\begin{equation}\label{eq10}
\lambda_j'= \int_{\gamma_j}\phi=
\int_{\sum a_{ij}\gamma_i}\phi=
\sum_{i=1}^{n}a_{ij}\lambda_i.
\end{equation}
Finally, let $j(\Lambda)={\Bbb Z}\lambda_1+\dots+{\Bbb Z}\lambda_{n}$
and $j'(\Lambda)={\Bbb Z}\lambda_1'+\dots+{\Bbb Z}\lambda_{n}'$.
Since $\lambda_j'= \sum_{i=1}^{n}a_{ij}\lambda_i$ and $(a_{ij})\in GL_{n}({\Bbb Z})$,
we conclude that $j(\Lambda)=j'(\Lambda)\subset {\Bbb R}$. In other words, the pseudo-lattices
$(\Lambda, {\Bbb R}, j)$ and $(\Lambda, {\Bbb R}, j')$ are isomorphic. Hence, 
$\alpha: {\cal V}\to {\cal PL}$  maps isomorphic complex algebraic curves to
the isomorphic pseudo-lattices, i.e. $\alpha$ is a functor.

Finally, let us show that $\alpha$ is a covariant functor. 
Indeed, let $\varphi_1,\varphi_2\in Mor ({\cal V})$. Then  
$\alpha(\varphi_1\varphi_2)= (\widetilde{\varphi_1\varphi_2})^*=\widetilde\varphi_1^*\widetilde\varphi_2^*=
\alpha(\varphi_1)\alpha(\varphi_2)$. Lemma \ref{lm5} follows.
$\square$

\bigskip\noindent
\underline{Category ${\cal PPL}$}. 
A  {\it projective pseudo-lattice} (of rank $n$) is a triple 
$(\Lambda, {\Bbb R}, j)$, where $\Lambda\cong {\Bbb Z}^n$ and $j: \Lambda\to {\Bbb R}$ is a
homomorphism. A morphism of projective pseudo-lattices
$(\Lambda, {\Bbb C}, j)\to (\Lambda, {\Bbb R}, j')$
is the identity
$j\circ\psi=\varphi\circ j'$,
where $\varphi$ is a group homomorphism and $\psi$ is an  ${\Bbb R}$-linear
map. (Notice that unlike the case of pseudo-lattices, $\psi$ is a scaling
map as opposite to an inclusion map. This allows to the two pseudo-lattices
to be projectively equivalent, while being distinct in the category ${\cal PL}_n$.) 
It is not hard to see that any isomorphism class of a projective pseudo-lattice 
contains a representative given by $j: {\Bbb Z}^n\to  {\Bbb R}$ such that
$j(1,0,\dots,0)=1,\quad
j(0,1,\dots,0)=\theta_1,\quad  \dots, \quad  j(0,0,\dots,1)=\theta_{n-1},$
  where $\theta_i$ are  positive reals.
The projective pseudo-lattices of rank $n$  make up a category, which we denote by ${\cal PPL}_n$.

\bigskip\noindent
\underline{Category $W$}. 
Finally, the toric $AF$-algebras ${\Bbb A}_{\theta}$,
modulo the stable isomorphism between them,  make up a category, 
which we shall denote by $W_n$. 
\begin{lem}\label{lm6}
Let $U_n\subseteq {\cal PPL}_n$ be a subcategory consisting of the projective
pseudo-lattices $PPL=PPL(1,\theta_1,\dots,\theta_{n-1})$ for which the Jacobi-Perron
fraction of the vector $(1,\theta_1,\dots,\theta_{n-1})$ converges to the vector.
Define a map  $\beta: U_n\to W_n$  by the formula
$PPL(1,\theta_1,\dots,\theta_{n-1})\mapsto {\Bbb A}_{\theta},$
where $\theta=(\theta_1,\dots,\theta_{n-1})$. 
Then $\beta$ is a bijective functor, which maps isomorphic projective pseudo-lattices 
to the stably isomorphic  toric $AF$-algebras. 
 \end{lem}
{\it Proof.} It is evident that $\beta$ is injective and surjective. Let
us show that $\beta$ is a functor. Indeed,  according to 
 \cite{E}, Corollary 4.7, every totally 
ordered abelian group of rank $n$ has form ${\Bbb Z}+\theta_1 {\Bbb Z}+\dots+ {\Bbb Z}\theta_{n-1}$.
The latter   is a projective pseudo-lattice $PPL$  from the category $U_n$. 
On the other hand, by  the Elliott theorem \cite{Ell1}, the $PPL$ defines a stable isomorphism class of 
the toric $AF$-algebra ${\Bbb A}_{\theta}\in W_n$. 
Therefore, $\beta$  maps isomorphic projective pseudo-lattices (from the set $U_n$) to the stably isomorphic toric 
$AF$-algebras,  and {\it vice versa}. Lemma \ref{lm6} follows.  
 $\square$

\medskip
Let $PL(\lambda_1,\lambda_2,\dots, \lambda_n)\in {\cal PL}_n$ and 
$PPL(1,\theta_1,\dots,\theta_{n-1})\in {\cal PPL}_n$. 
To finish the proof of theorem \ref{thm1}, it remains to prove the following lemma.
\begin{lem}\label{lm7}
Let $F: {\cal PL}_n\to {\cal PPL}_n$ be a map given by formula
$PL(\lambda_1,\lambda_2,\dots, \lambda_n)\mapsto PPL\left(1, {\lambda_2\over\lambda_1},\dots, {\lambda_n\over\lambda_1}\right).$
Then $Ker~F=(0,\infty)$ and  $F$ is a functor which maps  isomorphic pseudo-lattices to isomorphic projective 
pseudo-lattices. 
\end{lem}
{\it Proof.} Indeed, $F$ can be thought as   a map from ${\Bbb R}^n$ to ${\Bbb R}P^n$. Hence 
$Ker~F= \{\lambda_1 : \lambda_1>0\}\cong (0,\infty)$. The second part of lemma is evident.  
$\square$

\medskip
Setting $n=6g-6$ ($n=2$) for $g\ge 2$ ($g=1$) in the lemmas \ref{lm5}-\ref{lm7}, one gets the
conclusion of theorem \ref{thm1}.
$\square$

\bigskip\noindent
{\sf Acknowledgments.} 
I am  grateful to Yu. ~I. ~Manin for some helpful advices.

    

\vskip1cm

\textsc{The Fields Institute for Mathematical Sciences, Toronto, ON, Canada,  
E-mail:} {\sf igor.v.nikolaev@gmail.com}

\end{document}